\documentclass{article}%
\usepackage{amsfonts}
\usepackage{amssymb}
\usepackage{graphicx}
\usepackage{amsmath}%
\begin{document}

\begin{center}
{\huge Metrical Distortion,}

\ 

{\huge Exterior Differential}

\ 

{\huge and Gau\ss 's Lemma}

\ \ \ {\huge \ \ \ \ } \ \ \ \ \ \ \ \ \ 

\textbf{Stephan V\"{o}llinger} \ 

\bigskip\ \ \ \ 
\end{center}

\begin{quote}
\textbf{ stephan.voellinger@web.de}
\end{quote}

\ \ \ \ \ \ \ \ \ \ \ \ \ 

\begin{quote}
\textbf{Abstract}
\end{quote}

Gau\ss 's Lemma is revised by showing that the point set association of the
double tangential space with the tangential space of a Riemannian manifold is
not the identity. The latter point set association is called a metrical
distortion, an isometry that actually induces the geometry. The definition of
the exterior differential, which linearizes a mapping that points on the top
of a Riemannian manifold, is worked out concretely by covariant gradient
transport, what includes a differential slip in contrast to inner
differentials. A differential slip is a scalar gauge theory that considers the
reparametrization of different length scales. The metrical distortion is
determined by geodesically radial volume preservation, whereas the Riemannian
exponential mapping is determined by geodesically radial length preservation.
The theory is exemplified by the exterior geometry of the 2-sphere.\ \ 

\ \ 

\ 

\begin{quote}
\textbf{Content}
\end{quote}

\begin{quotation}
\textbf{1 Introduction}

\textbf{2 The Exterior Differential}

\textbf{3 Metrical Distortions}

\textbf{4 Exterior Volume}

\textbf{5} \textbf{Gau\ss 's Lemma}

\textbf{6 Exterior Geometry of the 2-Sphere}

\textbf{Bibliography}\ \ \ \ \ \ \ \ 
\end{quotation}

\ \ \ \ \ 

\begin{quote}
\textbf{Key words}
\end{quote}

General global diffeomorphism groups, metrical distortions, exterior
differential, Gau\ss 's lemma, differential slip, exterior volume.\ 

\begin{quote}
\textbf{Mathematics subject classification 2020}
\end{quote}

58A05 Differentiable manifolds, foundations.\newpage

\begin{quote}
\textbf{1 Introduction}
\end{quote}

What is the linearization in $p\in M$ of a mapping

\begin{center}
\hfill$\Theta:\mathbb{R}^{n}\rightarrow(M,\left\langle .,.\right\rangle ^{g}%
)$\hfill(1.1)\ 
\end{center}

that transports differentiable curves to differentiable curves and that points
on top of a $n$-dimensional Riemannian manifold [5] $(M,\left\langle
.,.\right\rangle ^{g})$?

As a reminder, the inner geometry of a Riemannian manifold, the
differentiability of chart transformations

\begin{center}
$f:=h\circ k^{-1}:\mathbb{R}^{n}|_{k(U\cap V)}\rightarrow\mathbb{R}^{n}$
\end{center}

of two different charts [5,6]

\begin{center}
$h:U\subset M\rightarrow\mathbb{R}^{n}$
\end{center}

and

\begin{center}
$k:V\subset M\rightarrow\mathbb{R}^{n}$
\end{center}

is treated by curve transportation. That is, if

\begin{center}
\hfill$v:=[t\rightarrow c(t)]$, $t\in\mathbb{R}$\hfill(1.2)

\hfill$c(t)\in M$, $c(0)=p$\hfill\ \ \ \ \ \ \ 
\end{center}

is the tangential vector equivalence class [5,6], alternatively defined as a
vector $n$-tuple

\begin{center}
\hfill$v^{i}:=\frac{d}{dt}|_{t=0}k^{i}\circ c(t)$, $v=v^{i}\partial_{i}$,
$i=1,..,n$,\hfill(1.3)

Sum convention! $v^{i}\partial_{i}=\underset{i}{%
{\textstyle\sum}
}v^{i}\partial_{i}$ \ \ \ \ \ \ \ \ \ \ 
\end{center}

in the image of $k$ at the outset, where $\partial_{i}$ denotes the coordinate
frame (see remark 1.1, I), then the inner differential $Df$ is defined by
inner curve transportation

\begin{center}
\hfill$Df[t\rightarrow k\circ c(t)]:=[t\rightarrow f\circ k\circ
c(t)]=[t\rightarrow h\circ c(t)]$\hfill(1.4)
\end{center}

what is $id$ on the manifold and

\begin{center}
\hfill$Dfv=(Df)_{i}^{j}v^{i}\partial_{j}=v^{j}\partial_{j}$\hfill(1.5)
\end{center}

in the image coordinate chart $h$, that is, the classical inner Jacobian in coordinates.

There are three methods to define tangential vectors in the inner geometry:
Equivalence classes of curves, coordinate $n$-tuples and derivations [6].
Derivations [5,6] are differential operators $\partial_{v}$ that operate on
germs of differentiable functions

\begin{center}
$\varphi_{p}:U_{p}\rightarrow\mathbb{R}$
\end{center}

defined in the vicinity $U_{p}$ of a point $p\in M$. In the inner geometry, if
$v$ is given as curve equivalence class (1.2), the associated derivation is

\begin{center}
\hfill$\partial_{v}\varphi_{p}:=\frac{d}{dt}|_{t=0}\varphi_{p}\circ
c(t)\in\mathbb{R}$.\hfill(1.6)
\end{center}

The result of a derivation is thereby a scalar function on the manifold. On
the other hand, if an abstract derivation $\partial_{v}$ is given at the
outset, it is associated to the vector $n$-tuple in the inner coordinate chart
$h$ by

\begin{center}
\hfill$(\partial_{v}h^{1},..,\partial_{v}h^{n})^{T}$\hfill(1.7)
\end{center}

what in turn defines a curve equivalence class [6]. The transposition
$(..)^{T}$ in (1.7) makes a row $n$-tuple to a column $n$-tuple.

\textbf{Remark 1.1 }(I) (Definition and covariance of coordinate frames.) A
global vector $v^{i}\partial_{i}$, where the index $i$ belongs to a coordinate
chart $h$, is invariant with respect to coordinate transformations. Thereby,
it follows necessarily that

\begin{center}
\hfill$\partial_{j}=(Df^{-1})_{j}^{i}\partial_{i}$,\hfill(1.8)
\end{center}

that is, the frame transforms covariant with the contragredient matrix
$(Df)^{-1,T}$. If a $n$-tuple transforms with a matrix, then the frame
transforms with the contragredient matrix in order to ensure the invariance.
With (1.8) a framed vector $n$-tuple is globally invariant

\begin{center}
\hfill$v^{i}\partial_{i}=\delta_{k}^{i}v^{k}\partial_{i}=(Df)_{k}^{j}%
(Df^{-1})_{j}^{i}v^{k}\partial_{i}=v^{j}\partial_{j}$.\hfill(1.9)
\end{center}

The frame is defined as the gradient equivalence class of the level set of the
$i$-th component function $h^{i}$ of the chart $h$

\begin{center}
\hfill$\partial_{i}:=[h^{i}]$,\hfill(1.10)
\end{center}

where the equivalence in $p\in M$

\begin{center}
$[\varphi_{p}^{1}]=[\varphi_{p}^{2}]$
\end{center}

of two function germs $\varphi_{p}^{1/2}:M\rightarrow\mathbb{R}$ holds by
definition exactly if the gradient in any chart $h$ is equal as covariant row
$n$-tuple

\begin{center}
$\partial_{i}(\varphi_{p}^{1}\circ h^{-1})=\partial_{i}(\varphi_{p}^{2}\circ
h^{-1})$, $i=1,..,n$.
\end{center}

(II) (Insufficiency of the inner differential.) Let (1.1) denote a mapping
that points on top of a Riemannian manifold, where the preimage space is
chosen without restriction to be trivial. If the exterior differential is
chosen to be the inner differential of the function $h\circ\Theta
:\mathbb{R}^{n}\rightarrow\mathbb{R}^{n}$, then a flat manifold is
automatically set in the background, namely the manifold with the inverse
Riemannian tensor

\begin{center}
$g^{ij}:=\underset{k}{%
{\textstyle\sum}
}(D(h\circ\Theta))_{k}^{i}(D(h\circ\Theta))_{k}^{j}$,

$g^{ik}g_{kj}=\delta_{j}^{i}$, $\left\langle v,w\right\rangle ^{g}%
:=g_{ij}v^{i}w^{j}$.
\end{center}

The latter definition does not linearize a mapping that points on top of a
curved manifold, explicit reparametrizations by the length parameter
associated to the line element are not considered. The differential acts on
frames which have the dual transformation behavior. The exterior concept of
covariance in definition 2.2 must naturally be gauged by a
reparametrization.\hfill\hfill\ \ \ \ \ \ \ \ \ \ $\blacksquare$\ 

Curve parameters must be well distinguished at the outset in different
manifolds and global differentiable mappings between them. The problem is now
how to fix the reparametrization concretely. This is done by the differential
slip of metrical distortions.

\textbf{Main Theorem 1.2 }Vectors transform covariant with respect to global mappings.

The metrical distortion $\Theta_{g}$ with origin $p$ of a Riemannian manifold
$(M,\left\langle .,.\right\rangle ^{g})$ is the point set association

\begin{center}
$\Theta_{g}|_{\Theta_{g}^{-1}(M\backslash cutlocus(p))}:T_{0}T_{p}M\rightarrow
M$
\end{center}

of the double tangential space with the manifold which is an isometry (theorem
3.1), where the domain of $\Theta_{g}^{-1}$ is restricted to the complement of
the cut locus with respect to $p$ [5].

The differential slip (definition 2.2, II) of the exterior differential
(definition 2.2, I) of $\Theta_{g}$ is determined by geodesically radial
volume preservation (theorem 5.2).

The covariant \ differential of any arbitrary mapping (1.1) can be decomposed
into a classical covariant differential $\mathbb{R}^{n}\rightarrow
T_{p}M\simeq\mathbb{R}^{n}$ without differential slip and the differential of
$\Theta_{g}$ with origin $p$ (lemma 2.4).\hfill\hfill
\ \ \ \ \ \ \ \ \ \ $\blacksquare$\ 

\begin{quote}
\textbf{2 The Exterior Differential}
\end{quote}

A tangential vector $v=v^{i}\partial_{i}\in T_{p}M$ can be seen as the level
set equivalence class of a function $\varphi_{p}$ in the vicinity of $p$. A
coordinate frame $\partial_{i}=[h^{i}]$ transforms covariant as in remark 1.1,
I.\ In analogy to define inner contravariant vectors as column $n$-tuples,
equivalence classes of curves and derivations as abstract differential
operators it is possible to define global covariant vectors which are framed
inner coordinate $n$-tuples.

Let $\exp_{p}:T_{p}M\rightarrow M$ denote Riemannian normal coordinates [5].
$\exp_{p}$ establishes an orthonormal contact coordinate system in the origin
$p$ for any $p\in M$.

\textbf{Proposition 2.1 }The following concepts of a tangential vector
$v=v^{i}\partial_{i}\in T_{p}M$ are equivalent

(I) (Level set equivalence class.) $v=[\varphi_{p}^{v}]$ is a level set
equivalence class with the level sets of the scalar function\ $\varphi_{p}%
^{v}:M\rightarrow\mathbb{R}$ in the vicinity of $p$.

(II) (Row $n$-tuple.) $v$ is represented by the global row $n$-tuple $(\xi
_{1},..,\xi_{n})$

\begin{center}
\hfill$D_{0}\exp_{p}^{-1}(v)=\underset{k}{%
{\textstyle\sum}
}\xi_{k}n_{k}\in\mathbb{R}_{synth}^{n}$\hfill(2.1)
\end{center}

with the synthetical coordinate frame $n_{k}$ and the synthetical space [1]

\begin{center}
$\mathbb{R}_{synth}^{n}:=T_{0}T_{p}M$.
\end{center}

Moreover, by definition

\begin{center}
\hfill$d\varphi_{p}^{v}:=(\partial_{1}(\varphi_{p}^{v}\circ\exp_{p}%
),..,\partial_{n}(\varphi_{p}^{v}\circ\exp_{p}))=(\xi_{1},..,\xi_{n})$.\hfill(2.2)
\end{center}

(III) (Covariant gradient.) $v$ is the covariant gradient

\begin{center}
\hfill$\varphi_{p}\rightarrow v(\varphi_{p}):=(\xi_{1}\partial_{1}\varphi
_{p},..,\xi_{n}\partial_{n}\varphi_{p})$,\hfill(2.3)

$\varphi_{p}:M\rightarrow\mathbb{R}$, \ \ \ \ \ 
\end{center}

that is, a componentwise scalar differential operator, where $\partial_{k}$,
$k\in\{1,..,n\}$ denotes the coordinate frame with respect to normal
coordinates $\exp_{p}^{-1}$.

\textbf{Proof }"(I)$\Leftrightarrow$(II)" Considering remark 1.1, I, the
reparametrization free natural identification [5]

\begin{center}
$D_{0}\exp_{p}:T_{0}T_{p}M\rightarrow T_{p}M$
\end{center}

needs to be stated for covariance. The freedom of the choice of the basis in
$T_{0}T_{p}M$ (synthetical frame [1]) makes it possible to choose a
reparametrization free initial basis. However, since it is an isometry the
musical isomorphisms $.^{\#}$ respectively $.^{\flat}$ [6] allow the dual
version. Concretely these isomorphisms are $id$ on the manifold and are
defined in coordinates by rising and lowering indexes via the metrical tensor

\begin{center}
$(v^{k})_{i}^{\flat}=g_{ik}v^{k}=v_{i}$, $g^{ik}\partial_{k}=d^{i}$
\end{center}

if $v^{\#}=v^{i}\partial_{i}$ respectively $v^{\flat}=v_{i}d^{i}$ denotes the
interior contravariant respectively covariant representation.

The covariant version of the natural identification is

\begin{center}
$D_{0}\exp_{p}[x^{k}]=[x^{k}\circ\exp_{p}^{-1}]$, $\forall k$
\end{center}

where $x^{k}:T_{p}M\rightarrow\mathbb{R}$ is the $k$-th component function
with respect to the choice of an orthonormal basis in the origin $p$. If
$\partial_{k}(p)\in$ $T_{p}M$ denotes the choice of the origin frame and
$v=v^{k}\partial_{k}$, then the component function in the vicinity of the
origin is defined by

\begin{center}
$x^{k}(v)=\left\langle \Pi_{\partial_{k}}v,\partial_{k}\right\rangle ^{g}$,
\end{center}

where $\Pi_{\partial_{k}}$ denotes the orthonormal projection on the
$1$-dimensional subspace. Moreover

\begin{center}
$[x^{k}]=n_{k}\in T_{0}T_{p}M$.
\end{center}

The differential is defined by transport of level sets as in definition 2.2
naturally without reparametrization.

The equivalence of I and II is established by (2.2) and

\begin{center}
$\varphi_{p}^{v}:=\underset{k}{%
{\textstyle\sum}
}\xi_{k}x^{k}\circ\exp_{p}^{-1}$.
\end{center}

If $\varkappa_{k}^{i}(p)$ denotes the Jacobian of the coordinate change
$h\circ\exp_{p}(0)$, $\forall p\in domain(h)\subset M$, then

\begin{center}
\hfill$v^{i}(h(p))=\underset{k}{%
{\textstyle\sum}
}\xi_{k}\varkappa_{k}^{i}(p)$.\hfill(2.4)
\end{center}

"(II)$\Leftrightarrow$(III)" The differential operator in III is given in an
arbitrary coordinate system by

\begin{center}
$\xi_{k}\varkappa_{k}^{i}\partial_{i}=\xi_{k}\varkappa_{k}^{j}(Df^{-1}%
)_{j}^{i}\partial_{i}=\xi_{k}\varkappa_{k}^{j}\partial_{j}$
\end{center}

with the coordinate change $f:i\rightarrow j$ and the normal coordinate index
$k$. The result of the latter differential operator on $\varphi_{p}$ in
coordinates is in any component scalar, what is elementarily computed with the
chain rule.\hfill\hfill\ \ \ \ \ \ \ \ \ \ $\blacksquare$\ 

A classical covariant transformation of a coordinate change

\begin{center}
$f:\mathbb{R}^{n}\rightarrow\mathbb{R}^{n}$, $f:i\rightarrow j$
\end{center}

is defined by

\begin{center}
\hfill$Df[\varphi_{p}^{v}]:=[\varphi_{p}^{v}\circ f^{-1}]$,\hfill(2.5)

$(Df^{-1})_{j}^{i}v_{i}=v_{j}$, $v=[\varphi_{p}^{v}]\in T_{p}\mathbb{R}^{n}$. \ \ \ 
\end{center}

The latter operation is the application of the contragredient matrix
$Df^{-1,T}$ on the covector $v_{i}$, $i=1,..,n$ to obtain the image vector
$v_{j}$, $j=1,..,n$.

The generalized differential for mappings that point on curved manifolds
includes a reparametrization in order to consider the different length parameters.

\textbf{Definition 2.2 }(I) (Exterior differential.) Let (1.1) denote a
bijective mapping $\Theta(x)=p$ that points on top of a Riemannian manifold.
The \emph{exterior differential }on the tangential space

\begin{center}
$D\Theta:T_{x}\mathbb{R}^{n}\rightarrow T_{p}M$
\end{center}

is defined by gradient transport

\begin{center}
\hfill$D\Theta\lbrack\varphi_{x}]:=[\varphi_{x}\circ\Theta^{-1}]$, \hfill(2.6)

$\varphi_{x}:\mathbb{R}^{n}\rightarrow\mathbb{R}$, in a vicinity of $x$. \ \ \ \ 
\end{center}

With the notion in proposition 2.1 in the preimage space

\begin{center}
$[\varphi_{x}]=d\varphi_{x}=(\partial_{1}\varphi_{x},..,\partial_{n}%
\varphi_{x})$,

$\partial_{k}\varphi_{q}:=\frac{d}{ds}|_{s=0}\varphi_{q}(q+se_{k})$, $e_{k}$
unit vector
\end{center}

and in the image space

\begin{center}
$[\varphi_{x}\circ\Theta^{-1}]=d(\varphi_{x}\circ\Theta^{-1})$

$:=(\partial_{1}(\varphi_{x}\circ\Theta^{-1}\circ\exp_{p}),..,\partial
_{n}(\varphi_{x}\circ\Theta^{-1}\circ\exp_{p}))$,

$\partial_{i}(\varphi_{x}\circ\Theta^{-1}\circ\exp_{p})(p):=\frac{d}%
{dt}|_{t=0}\varphi_{x}\circ\Theta^{-1}\circ\exp_{p}(\exp_{p}^{-1}(p)+te_{i})$.
\end{center}

$s$ is the length parameter in the preimage space and $t$ is the length
parameter on the manifold.

(II) (Differential slip.) The curve parameters $s$ in the preimage space and
$t$ in the image space are different. The differential $\frac{dt}{ds}$ along a
curve is called \emph{differential slip} along the curve. The differential
slip is the isometry that induces the base point change

\begin{center}
$T_{0}T_{p}M\simeq T_{0}T_{q}M$
\end{center}

as in lemma 2.3. It is the gauge factor for the choice of an origin in
$\mathbb{R}_{synth}^{n}$ (lemma 2.4).

(III) (Time change.) Let $M$ be a manifold and $c_{M}:[0,t_{\max}]\rightarrow
M$ a curve parametrized by $t$. A \emph{differentiable time change}
$t\rightarrow s$ is a positive, continuous scalar function

\begin{center}
$a:M\rightarrow\mathbb{R}$, $a>0$,
\end{center}

such that

\begin{center}
\hfill$t(s):=\overset{s}{\underset{0}{%
{\textstyle\int}
}}a(c_{M}(l))dl$,\hfill(2.7)
\end{center}

for arbitrary curves through the point of a manifold. It is a scalar gauge
theory independent from the direction of the curve.\hfill\hfill
\ \ \ \ \ \ \ \ \ \ $\blacksquare$\ 

The properties of a differential are preserved by an additional time change
structure. By definition of the exterior differential, the differential slip
is a time change.

\textbf{Lemma 2.3 }The base point change

\begin{center}
\hfill$D_{q}\exp_{p}:(T_{q}T_{p}M,\left\langle .,.\right\rangle _{T_{q}T_{p}%
M}^{eukl})\rightarrow(T_{q}M,\left\langle .,.\right\rangle _{T_{q}M}^{g}%
)$\hfill(2.8)

$\simeq(T_{0}T_{q}M,\left\langle .,.\right\rangle _{T_{0}T_{q}M}^{eukl})$
\end{center}

is a linear coordinate change isometry.

\textbf{Proof }Let $v\in T_{p}M$ and$\ \exp_{p}(v)=q\in M$. As in [1, lemma
3], the exterior differential of the mapping

\begin{center}
$\tau_{v}:T_{p}M\rightarrow T_{p}M$,

$m\rightarrow m+v$
\end{center}

serves as the parallel translation $id_{pq}$ in the affine space $TT_{p}M$ [1]

\begin{center}
$D\tau_{v}=id_{pq}:T_{p}T_{p}M\rightarrow T_{q}T_{p}M$.
\end{center}

The origin in $T_{q}T_{p}M$ is the point that is associated with $q$ in the
manifold and

\begin{center}
$D_{v}\Theta_{q}^{p}=D_{v}\exp_{p}\circ id_{pq}:T_{p}T_{p}M\rightarrow
T_{q}T_{p}M\rightarrow T_{q}M$
\end{center}

with

\begin{center}
$\Theta_{q}^{p}:w\rightarrow\exp_{p}(v+w)$.
\end{center}

Follows because $\exp_{q}^{-1}\circ\exp_{p}$ is a coordinate change which is
$id$ on the manifold. See [1, lemma 3 III] for details.\hfill\hfill
\ \ \ \ \ \ \ \ \ \ $\blacksquare$\ 

With lemma 2.3 $TT_{p}M$ builds a natural affine space over the manifold [1].

\textbf{Lemma 2.4 }(I) (Reparametrized gradient transport of $D_{q}\exp_{p}$.)
Let $s$ denote the length parameter in the affine space $TT_{p}M$ and $t$
denote the Riemannian length on the manifold $M$.

Let $\varphi:T_{p}M\rightarrow\mathbb{R}$ denote a scalar function on the
tangential space in the vicinity of $q=\exp_{p}(v)$, such that $[\varphi]\in
T_{q}T_{p}M$ and

\begin{center}
$\varphi^{\prime}((\exp_{p}^{-1})^{1},..,(\exp_{p}^{-1})^{n}):=\varphi
(\exp_{p}^{-1})$
\end{center}

the associated function in the interior normal coordinate system such that

\begin{center}
$[\varphi^{\prime}((\exp_{p}^{-1})^{1},..,(\exp_{p}^{-1})^{n})]=[\varphi
^{\prime}]\in T_{q}M$.
\end{center}

Then, $D_{q}\exp_{p}$ is the isometry as in (2.8) exactly if it is computed by
reparametrized gradient transport

\begin{center}
\hfill\ $D_{q}\exp_{p}[\varphi]\overset{\text{definition 2.2, I}}{=}%
[\varphi\circ\exp_{p}^{-1}]\overset{\text{Isometry}}{=}\frac{dt}{ds}%
[\varphi^{\prime}]$.\hfill(2.9)
\end{center}

(II) (Renormalized curve transportation of $D_{q}\exp_{p}$.) Let $c(t)\in M$
denote a curve with $c(0)=q$. Reparametrized gradient transport is equivalent
to renormalized curve transportation

\begin{center}
\hfill\ $D_{q}\exp_{p}[s\rightarrow\exp_{p}^{-1}\circ c(t(s))]=[t\rightarrow
c(t)]$,\hfill\ \ \ \ \ \ \ \ \ \ 

$c(t)\in M$, $c(0)=q$,
\end{center}

where, since $\exp_{p}$ is an isometry, it commutes with the musical
isomorphisms $.^{\#}$ respectively $.^{\flat}$ as in the proof of proposition
2.1, such that the application of the diffeomorphism on a curve equivalence
class can be defined by

\begin{center}
$D\exp_{p}v^{\#}:=(D\exp_{p}v)^{\#}$.
\end{center}

(III) (Decomposition of an arbitrary covariant differential in a composition
of a classical differential and $D_{q}\exp_{p}$.) If $\Theta$ is a
diffeomorphism that points on top of a manifold as in (1.1), then $\exp
_{p}^{-1}\circ\Theta$ is a classical diffeomorphism $\mathbb{R}^{n}\rightarrow
T_{p}M\simeq\mathbb{R}^{n}$ and

\begin{center}
\hfill$D_{q}\Theta=D_{q}\exp_{p}\circ D_{q}(\exp_{p}^{-1}\circ\Theta)$,\hfill(2.10)
\end{center}

such that the exterior differential of $\Theta$ can be computed by the
computation of a classical covariant transformation (2.5) and the computation
of the differential of the isometry $\exp_{p}$ as in I.

\textbf{Proof }(I) According to lemma 2.3 the by $D_{q}\exp_{p}$ induced basis
is orthonormal and if this basis is chosen as origin for $T_{q}M\simeq
\mathbb{R}^{n}$, $D_{q}\exp_{p}$ is $id$ on the level of covariant row tuples.

If $[\varphi]\in T_{q}T_{p}M$, with $domain(\varphi)\in T_{p}M$, then
according to proposition 2.1 with $TT_{p}M$ as manifold, $[\varphi]$ as
preimage row tuple is computed by

\begin{center}
$(\xi_{1},..,\xi_{n})_{i}=\frac{d}{ds}|_{s=0}\varphi\circ\exp_{p}^{-1}%
\circ\exp_{q}(\exp_{q}^{-1}(q)+se_{i})$,
\end{center}

since lemma 2.3 is exactly an isometry if the normal coordinate transformation
$\exp_{q}^{-1}\circ\exp_{p}$ is an orthonormal contact transformation to the
origin of $T_{q}T_{p}M$.

The row $n$-tuple associated to $[\varphi^{\prime}]$ "on top of the manifold"
$(M,\left\langle .,.\right\rangle ^{g})$ is computed by definition in the
normal coordinate system with origin $q$ by

\begin{center}
$(d(\varphi^{\prime}((\exp_{p}^{-1})^{1},..,(\exp_{p}^{-1})^{n})))_{i}%
=\frac{d}{dt}|_{t=0}\varphi\circ\exp_{p}^{-1}\circ\exp_{q}(\exp_{q}%
^{-1}(q)+te_{i}))$.
\end{center}

By time change $\frac{dt}{ds}\frac{d}{dt}=\frac{d}{ds}$ the latter row
$n$-tuple is exactly equal to $(\xi_{1},..,\xi_{n})$ if it is multiplied with
$\frac{dt}{ds}$, that is, if (2.9) holds.

(II) The component function $x_{q}^{i}$ on $T_{p}M$ in the vicinity of $q$ is
defined by

\begin{center}
$x_{q}^{i}(v+w)=\left\langle \Pi_{\partial_{i}}w,\partial_{i}\right\rangle
_{T_{p}M}^{g}$,

$\Pi_{\partial_{i}}$ orthonormal projection on $\partial_{i}$
\end{center}

if $\partial_{i}(p)\in$ $T_{p}M$ denotes the choice of an origin frame and
$v=v^{i}\partial_{i}$, $w=w^{i}\partial_{i}$ respectively $\exp_{p}(v)=q$.

If $v\in T_{q}M$ with

\begin{center}
$v^{\#}=(v^{1},..,v^{n})^{T}=[t\rightarrow c(t)]$
\end{center}

and the component function $[x_{q}^{i}]\in T_{q}T_{p}M$, with $domain(x_{q}%
^{i})\in T_{p}M$, then

\begin{center}
$\left\langle D_{q}\exp_{p}[x_{q}^{i}],v^{\#}\right\rangle _{T_{q}M}%
^{g}\overset{\text{isometry}}{=}\left\langle [x_{q}^{i}],D_{q}\exp_{p}%
^{-1}v^{\#}\right\rangle _{T_{q}T_{p}M}^{eukl}$

$=(D_{q}\exp_{p}^{-1}v^{\#})^{i}\in\mathbb{R}$, $\forall i$
\end{center}

and according to (2.9)

\begin{center}
$\left\langle D_{q}\exp_{p}[x_{q}^{i}],v^{\#}\right\rangle _{T_{q}M}%
^{g}=\left\langle \frac{dt}{ds}\partial_{i},[t\rightarrow c(t)]\right\rangle
_{T_{q}M}^{g}$

$=\frac{dt}{ds}[t\rightarrow(\exp_{p}^{-1}\circ c(t))^{i}]\in\mathbb{R}$,
$\forall i$.
\end{center}

It can be concluded that

\begin{center}
$D_{q}\exp_{p}^{-1}v^{\#}=\frac{dt}{ds}[t\rightarrow\exp_{p}^{-1}\circ
c(t)]\in T_{q}T_{p}M$.
\end{center}

But on the level of contravariant column $n$-tuples

\begin{center}
$\frac{dt}{ds}\frac{d}{dt}\exp_{p}^{-1}\circ c(t)=\frac{d}{ds}\exp_{p}%
^{-1}\circ c(t(s))$,
\end{center}

such that

\begin{center}
$\frac{d}{ds}|_{s=0}\exp_{p}^{-1}\circ c(t(s))$
\end{center}

is the image of $D_{q}\exp_{p}^{-1}v^{\#}$, as column $n$-tuple in $T_{q}%
T_{p}M$ (which is equivalently a row $n$-tuple since the metric is trivial).

(III) The mapping

\begin{center}
$\exp_{p}^{-1}\circ\Theta:(\mathbb{R}^{n},\left\langle .,.\right\rangle
^{Eucl})\rightarrow(T_{p}M,\left\langle .,.\right\rangle _{T_{p}M}^{g})$
\end{center}

is a classical covariant mapping $\mathbb{R}^{n}\rightarrow T_{p}%
M\simeq\mathbb{R}^{n}$ since an orthonormal basis was chosen in $T_{p}M$ as
origin. The differential slip is trivial. The classical covariant
diffeomorphism is computed as in (2.5) and $D_{q}\exp_{p}$ as in I by
reparametrized gradient transport or by renormalized curve transportation
respectively. The covariant differential obeys the chain rule.\hfill
\hfill\ \ \ \ \ \ \ \ \ \ $\blacksquare$\ 

\begin{quote}
\textbf{3 Metrical Distortions}
\end{quote}

Let $\varkappa(q)$, $q\in M$ be defined as in (2.4) in an arbitrary coordinate
chart $h$, that is, by the inner differential $Dh\circ\exp_{q}(0)$ with the
Riemann exponential mapping with origin $q$, $\forall q\in domain(h)\subset
M$. Since the metric is trivial in the origin of normal coordinates it follows that

\begin{center}
$g^{ij}:=\underset{k}{%
{\textstyle\sum}
}\varkappa_{k}^{i}\varkappa_{k}^{j}$, $g^{ik}g_{kj}=\delta_{j}^{i}$,
\end{center}

for the metric in the chart $h$.

It is now possible to pose the problem to find a solution of the equation

\begin{center}
\hfill$D\Theta(x)\partial_{k}=\varkappa_{k}^{j}(\Theta(x))\partial_{j}$,\hfill(3.1)
\end{center}

if the Riemannian manifold $(M,\left\langle .,.\right\rangle ^{g})$ is
prescribed as data of the PDE. That is, to find a point set association
$\Theta$ that actually induces the prescribed geometry.

The problem was solved in [1] to be integrated by the natural point set
association of the double tangential space with origin $p\in M$ with the manifold

\begin{center}
\hfill$\Theta_{g}|_{\Theta_{g}^{-1}(M\backslash cutlocus(p))}:T_{0}%
T_{p}M\rightarrow M$.\hfill(3.2)
\end{center}

The differential of $\exp_{p}$ is eventually singular at the cut locus,
thereby the image is restricted to the complement set of $cutlocus(p)\subset
M$ [5] which is a Lebesgue zero set [5].

\textbf{Theorem 3.1 }(I) (Isometry property.) By construction, the Levy-Civita
metrical distortion $\Theta_{g}$ is an exterior diffeomorphism and solves
(3.1), the geometric PDE [1, theorem 14]. (3.1) is the PDE, which actually
induces the geometry since it is an isometry

\begin{center}
\hfill$\left\langle D\Theta_{g}(x)v,D\Theta_{g}(x)w\right\rangle
^{g}=\left\langle v,w\right\rangle ^{T_{q}T_{p}M}$,\hfill(3.3)

$v,w\in T_{q}T_{p}M$ \ \ \ \ \ 
\end{center}

by construction, with the trivial Euclidean metric in the affine space
$TT_{p}M$ with origin $x$ if $\Theta_{g}(x)=q$.

(II) (Geodesically radial mapping.) Let $exp_{p}$ denote the Riemannian
exponential function [5]. The point set association (3.2) is a geodesic
mapping, that is, Euclidean geodesics with origin $0$ are transferred to
geodesics in the manifold with origin $p$

\begin{center}
\hfill$\Theta_{g}(exp_{0}^{\mathbb{R}_{synth}^{n}}(rv))=exp_{p}(r^{\prime
}(r,v)v)$\hfill(3.4)
\end{center}

where

\begin{center}
$v\in T_{p}M\simeq$ $\mathbb{R}_{synth}^{n}$,$\left\Vert v\right\Vert =1$ \ \ \ \ 
\end{center}

is a unit vector $n$-tuple in the origin and $r^{\prime}(r,v)$ is the radial
contraction along the geodesic with initial vector $v$.

$\Theta_{g}$ is determined by a point set association

\begin{center}
\hfill$(r,v)\rightarrow(r^{\prime}(r,v),v)$\hfill(3.5)
\end{center}

of spherical polar coordinates in synthetical coordinates to spherical polar
coordinates of normal coordinates in the manifold, where the vector

\begin{center}
$v=v(\varphi_{1},..,\varphi_{n-1})$
\end{center}

is the preimage of spherical polar coordinates with the spherical angles
$\varphi_{i}$, $i=1,..,n-1$ and radius $1$.

The mapping in arbitrary coordinates is called geodesically radial with origin
$p$.

(III) (Equivalence of the radial contraction and the differential slip.) The
radial contraction $\frac{dr^{\prime}}{dr}$ is equal to the differential slip
$\frac{dt}{ds}$.

\textbf{Proof }(I) Let $w\in T_{p}M$ and$\ \exp_{p}(w)=q\in M$. As in [1,
lemma 3], the exterior differential of the mapping

\begin{center}
$\tau_{w}:T_{p}M\rightarrow T_{p}M$,

$m\rightarrow m+w$
\end{center}

serves as the parallel translation $id_{pq}$ in the affine space $TT_{p}M$

\begin{center}
$D\tau_{w}=id_{pq}:T_{p}T_{p}M\rightarrow T_{q}T_{p}M$,
\end{center}

where the origin $x$ in $T_{q}T_{p}M$ is the point that is associated by
$\Theta_{g}(x)=q$ in the manifold.

The natural identification of tangential vectors in the origin of normal
coordinates [5]

\begin{center}
$D_{0}\exp_{p}:\mathbb{R}_{synth}^{n}=T_{0}T_{p}M=T_{p}T_{p}M\rightarrow
T_{p}M$
\end{center}

can be stated covariant as in proposition 2.1 for any $p\in M$.

The exterior differential of $\Theta_{g}$ naturally resembles the differential
of the Riemannian exponential function which is an isometry as in lemma 2.3
respectively [1, lemma 3] for more details.

(II) The set of points in $TT_{p}M$ defines an affine space as in [1], the
affine space over the manifold. The choice of coordinates in $TT_{p}M$ is
synthetically induced. According to lemma 2.3 the affine space over the
manifold is orthonormal to the manifold. By parallel translation of tangential
vectors it follows that the curve that is associated with a straight line from
$p$ in $T_{p}M$ is a straight line in the affine space and that the initial
angle of $\frac{w}{\left\Vert w\right\Vert _{\mathbb{R}^{n}}}=v$ with respect
to the orthonormal frame of normal coordinates is the same in the affine space.

(III) Since the differential slip is a time change, it is independent from the
direction of a curve through any point. A radial curve through the point
exists and has the contraction $\frac{dr^{\prime}}{dr}$.\hfill\hfill
\ \ \ \ \ \ \ \ \ \ $\blacksquare$

For a complete determination of $\Theta_{g}$ it is thus sufficient to state
$r^{\prime}(r,v)$, which determines the radial contraction. That will be
accomplished by volume preservation.

\textbf{Remark 3.2 }(Turbulences [1,2,3].) If $\varkappa$ is associated to the
Levy-Civita geodesic system, any $\varkappa\circ O$ with an orthonormal matrix
$O$ such that the matrix multiplication

\begin{center}
\hfill$(g^{ij})_{ij}=(\underset{k}{%
{\textstyle\sum}
}\varkappa_{k}^{i}\varkappa_{k}^{j})_{ij}=\varkappa\circ\varkappa
^{T}=(\varkappa\circ O)\circ(\varkappa\circ O)^{T}$\hfill(3.6)
\end{center}

is preserved, is associated to autoparallel geodesics and autoparallel
exponential functions.

The theory of the geometric PDE (3.1) can be generalized to arbitrary
$\varkappa\circ O$ such that (3.6) is obeyed and autoparallel metrical
distortions $\Theta_{\varkappa\circ O}$ [1,2,3] are obtained.

The mapping

\begin{center}
\hfill$\Theta_{\varkappa\circ O}^{-1}\circ\Theta_{g}:\mathbb{R}_{synth}%
^{n}\rightarrow\mathbb{R}_{O-synth}^{n}$\hfill(3.7)
\end{center}

is called a \emph{turbulence}.

A turbulence is a classical covariant mapping which has an orthonormal
differential.\hfill\hfill\ \ \ \ \ \ \ \ \ \ $\blacksquare$

\begin{quote}
\textbf{4 Exterior Volume}
\end{quote}

If $(M,\left\langle .,.\right\rangle ^{g})$ is a $n$-dimensional Riemannian
manifold and $h:U\rightarrow\mathbb{R}^{n}$, $U\subset M$ a chart, then

\begin{center}
\hfill$Vol_{int}(U):=\underset{h(U)}{%
{\textstyle\int}
}dvol_{int}^{g}(z)=\underset{h(U)}{%
{\textstyle\int}
}\det(g_{ij})^{ij}dz$\hfill(4.1)
\end{center}

is called the \emph{interior volume}. The interior volume is independent of
the chart [5] and defines thereby a measure on the manifold with respect to
the Borelian measure space [7].

\textbf{Definition 4.1 }(Exterior Volume.) Let $(M,\left\langle
.,.\right\rangle ^{g})$ denote a Riemannian manifold of dimension $n$. Choose
normal coordinates with origin $q\in M$ parametrized by $z$. Let

\begin{center}
$\varphi^{\prime}(z)\in C_{0}^{0}(\mathbb{R}^{n})$ (continuous compact support),

$\int\varphi^{\prime}(z)dz=1$
\end{center}

and

\begin{center}
$\varphi_{\varepsilon}^{\prime}(z):=\frac{1}{\varepsilon^{n}}\varphi^{\prime
}(\frac{z}{\varepsilon})$
\end{center}

such that

\begin{center}
$\int\varphi_{\varepsilon}^{\prime}(z)dz=1$
\end{center}

and

\begin{center}
$\varphi_{\varepsilon}^{\prime}\underset{\varepsilon\rightarrow0}{\rightarrow
}\delta_{0}^{\prime}$,
\end{center}

with the delta distribution in $0$. Then

\begin{center}
$\varphi_{\varepsilon,q}:=\varphi_{\varepsilon}^{\prime}\circ\exp_{q}^{-1}$
\end{center}

converges with $\varepsilon\rightarrow0$ to the $\delta\emph{-}$%
\emph{distribution on the manifold}.

Any measure $vol_{ext}^{g}$ on the Borelian measure space [7] with respect to
the abstract point set $M$, which is absolutely continuous [7] with respect to
the interior measure, is called \emph{exterior volume} if for any $q\in M$

\begin{center}
\hfill$\underset{\varepsilon\rightarrow0}{\lim}\underset{M}{%
{\textstyle\int}
}\varphi_{\varepsilon,q}dvol_{ext}^{g}=1$.\hfill(4.2)
\end{center}

Any mapping $\Theta:\mathbb{R}^{n}\rightarrow M$ such that the image measure
of the Lebesgue measure obeys (4.2) is called \emph{infinitesimally volume
preserving}.\hfill\hfill\ \ \ \ \ \ \ \ \ \ $\blacksquare$

For any point $q\in M$, the exterior volume is infinitesimally equal to the
Lebesgue volume of the associated set in normal coordinates with origin $q$.

\textbf{Remark 4.2 }(Exterior Volume and interior volume.) There is no reason
why exterior and interior volume should be equivalent. If there is an
infinitesimal exhaustion which is linearizable by parallelepipeds on top of
the manifold, there is no reason why the tesselation by parallelepipeds should
also fit infinitesimally in the interior exhaustion with respect to a
coordinate chart. The linearization differs by a differential slip which is a
scalar gauge theory of holonomic coordinate grids and does not fit in the
interior exhaustion.\hfill\hfill\ \ \ \ \ \ \ \ \ \ $\blacksquare$

The exterior volume is a restatement of the volume of a Riemannian manifold.

\textbf{Theorem 4.3 }(Determination of the Exterior Volume.) Any exterior
measure according to definition 4.1 is uniquely determined by the image
measure of the metrical distortion $\Theta_{g}$ which is thereby
infinitesimally volume preserving. For any connected open set $O\subset
\mathbb{R}_{synth}^{n}$ within the complement set of the cut locus with
respect to the origin of $\Theta_{g}$, the transformation formula

\begin{center}
\hfill$\underset{O}{\int}dvol^{eucl}=\underset{\Theta_{g}(O)}{\int}%
dvol_{ext}^{g}$\hfill(4.3)
\end{center}

holds.

\textbf{Proof }Let $(\Theta_{g})_{\#}\lambda$ denote the image measure of the
Lebesgue measure in synthetical coordinates such that $\Theta_{g}(0)=q\in M$.
$(\Theta_{g})_{\#}\lambda$ is absolutely continuous with respect to the
interior volume [7]. Let

\begin{center}
$\varphi^{synth}\in C_{0}^{0}(\mathbb{R}_{synth}^{n})$ (continuous compact support),

$\int\varphi^{synth}(x)dz=1$
\end{center}

and

\begin{center}
$\varphi_{\varepsilon}^{synth}(x):=\frac{1}{\varepsilon^{n}}\varphi
^{synth}(\frac{x}{\varepsilon})$
\end{center}

such that

\begin{center}
$\int\varphi_{\varepsilon}^{synth}(z)dz=1$.
\end{center}

Then, with the $\delta\emph{-}$distribution in synthetical coordinates and the
image measure theorem

\begin{center}
$1=\underset{\varepsilon\rightarrow0}{\lim}\underset{\mathbb{R}_{synth}^{n}}{%
{\textstyle\int}
}\varphi_{\varepsilon}^{synth}(\Theta_{g}^{-1}(y))d(\Theta_{g})_{\#}%
\lambda(y)$,
\end{center}

such that, since $q$ is arbitrary, it suffices to prove that $\varphi
_{\varepsilon}^{synth}\circ\Theta_{g}^{-1}$ converges to the $\delta\emph{-}%
$distribution in $q$ on the manifold according to definition 4.1. With the
geodesic deviation $\overrightarrow{\Theta}_{g}^{-1}$ in $q$

\begin{center}
$\Theta_{g}^{-1}\circ\exp_{q}(y)=:id(y)+\overrightarrow{\Theta}_{g}^{-1}(y)$,
\end{center}

where the natural identification $T_{0}T_{q}M=T_{q}M$ is used. It follows that

\begin{center}
$\underset{\varepsilon\rightarrow0}{\lim}\int\varphi_{\varepsilon}%
^{synth}\circ\Theta_{g}^{-1}\circ\exp_{q}(y)f(y)dy$

$=\underset{\varepsilon\rightarrow0}{\lim}\int\varphi_{\varepsilon}%
^{synth}(y+\overrightarrow{\Theta}_{g}^{-1}(y))f(y)d(y)$

$\underset{y=\varepsilon z}{\overset{\text{Substitution}}{=}}\underset
{\varepsilon\rightarrow0}{\lim}\int\varphi^{synth}(z+O(\varepsilon
))f(\varepsilon z)dz$

$\overset{\text{Maj.}}{\underset{\text{convergence [5]}}{=}}f(0)$, $\forall
f\in C_{0}^{0}$,
\end{center}

what establishes the convergence to the $\delta\emph{-}$distribution in $q$.
The Landau symbol

\begin{center}
$\underset{\varepsilon\rightarrow0}{\lim}O(\varepsilon)=0$
\end{center}

is introduced, because the derivative of the geodesic deviation is computed by

\begin{center}
$\underset{\varepsilon\rightarrow0}{\lim}\frac{1}{\varepsilon}\overrightarrow
{\Theta}_{g}^{-1}(\varepsilon z)=0$.
\end{center}

A geodesic is tangential in the origin of normal coordinates.\hfill
\hfill\ \ \ \ \ \ \ \ \ \ $\blacksquare$

The exterior volume is "bootstrapped" by the metrical distortion.

\textbf{Example 4.4 }(Radial mappings.) Let $(M,\left\langle .,.\right\rangle
^{g})$ denote a manifold. A geodesically radial mapping

\begin{center}
$f_{r^{\prime}}(x)=r^{\prime}(\left\vert x\right\vert )\frac{x}{\left\vert
x\right\vert }$,
\end{center}

with image in normal coordinates and

\begin{center}
$\frac{dr^{\prime}(s)}{ds}=a$
\end{center}

has the image classically differentiable volume element

\begin{center}
$\det(Df_{r^{\prime}}^{-1})=\frac{1}{(\frac{r^{\prime}(\left\vert x\right\vert
)}{\left\vert x\right\vert })^{n-1}a}$.
\end{center}

The global integral equation for radial mappings

\begin{center}
$%
{\textstyle\int}
\varphi_{p}\circ f_{r^{\prime}}(x)dx=%
{\textstyle\int}
\varphi_{p}(z)dvol_{int}^{g}(z)$, $\forall\varphi_{p}$, $p\in M$
\end{center}

is only solvable for dimension $2$ if $a=\frac{\left\vert x\right\vert
}{r^{\prime}g}$, as in theorem 5.2, where $g$ is the $2$-dimensional
restriction volume element in normal coordinates.

The exterior volume is equivalent to the image volume of $f_{r^{\prime}}$, if
the image of $f_{r^{\prime}}$ is identified with the abstract point set $M$.
The classical image volume of $f_{r^{\prime}}$ in normal coordinates is
equivalent to the exterior volume, where equivalence means equivalence up to a
scalar gauge factor due to a differential slip. The volumemetrical
approximation is however different since the global grid on top of the
manifold is different.\hfill\hfill\ \ \ \ \ \ \ \ \ \ $\blacksquare$

\begin{quote}
\textbf{5 Gau\ss 's Lemma}
\end{quote}

The classical Gau\ss 's lemma [5] can be stated as in the following theorem.

\textbf{Theorem 5.1 }(Classical Gau\ss 's lemma [5].) Let $exp_{p}$ denote the
Riemannian exponential function [5], with the choice of an orthonormal basis
in the origin such that $T_{p}M\simeq\mathbb{R}^{n}$. Let $q(t):=\exp_{p}%
(tv)$, $\left\Vert v\right\Vert =1$. Let $h$ denote a chart in the Riemannian
manifold $(M,\left\langle .,.\right\rangle ^{g})$ such that

\begin{center}
$h\circ$ $\exp_{p}:\mathbb{R}^{n}\rightarrow\mathbb{R}^{n}$
\end{center}

is a classically differentiable coordinate change diffeomorphism in the
vicinity of $p$.

Then, in the vicinity of $p$

\begin{center}
\hfill$D_{tv}h\circ\exp_{p}:(T_{tv}\mathbb{R}^{n},\left\langle
.,.\right\rangle ^{g})\rightarrow(T_{h(q(t))}\mathbb{R}^{n},\left\langle
.,.\right\rangle ^{g})$,\hfill(5.1)
\end{center}

is a linear isometry between the metrical spaces in coordinates equipped with
the metric $g_{ij}$, where $\exp_{p}(tv)=q(t)$ is the arc length parametrized
geodesic with initial unit tangential vector $v$. The parameter $t$ is the
Riemannian length parameter independent from the choice of coordinates. Especially

\begin{center}
\hfill$D_{tv}h\circ\exp_{p}(v)=\frac{d}{dt}h(q(t))$.\hfill(5.2)
\end{center}

$h\circ$ $\exp_{p}$ is the geodesically radial length preserving mapping,
where radial refers to the choice of the origin $p$.\hfill\hfill
\ \ \ \ \ \ \ \ \ \ $\blacksquare$

The coordinate change between normal coordinates must be computed by solving
the geodesic problems in the interior geometry. The latter problem must be
distinguished from exterior differentiation of the mapping

\begin{center}
$\exp_{p}:T_{p}\mathbb{R}^{n}\rightarrow(M,\left\langle .,.\right\rangle
^{g})$
\end{center}

which includes a differential slip.

As in chapter 4, the measure and integration theory on top of the manifold
$(M,\left\langle .,.\right\rangle ^{g})$

\begin{center}
$\underset{M}{%
{\textstyle\int}
}f(m)dvol^{g}(m)$, $m\in M$, $f:M\rightarrow\mathbb{R}$ measurable
\end{center}

is defined locally linearized for parallelepipeds in normal coordinates with
local origin as trivial volume element $dvol^{g}$ and then globally summarized.

The uniqueness of a geodesically radial volume preserving mapping determines
the generalized Gau\ss 's\ lemma.

\textbf{Theorem 5.2 }(Generalized Gau\ss 's\ lemma.) In the vicinity of $p$,
or more precisely in the complement set of the cut locus with respect to
$p$,\ the canonical point set association

\begin{center}
\hfill$\Theta_{g}:T_{0}T_{p}M\rightarrow M$,\hfill(5.3)
\end{center}

the Levy-Civita metrical distortion (without turbulence [1,2,3]), is the
geodesically radial infinitesimally volume preserving mapping (definition
4.1), which is an isometry as in (3.3).

The isometry property (theorem 3.1, I) is $2$-dimensional restrictable. In the
case of dimension $n=2$ the metrical distortion is globally volume preserving,
it transforms the Lebesgue measure in synthetical coordinates to the inner
volume of the Riemannian manifold.

\textbf{Proof }$\Theta_{g}$ is geodesically radial by theorem 3.1. By the
isometry property (3.3) the metrical distortion is infinitesimally volume
preserving (theorem 4.3).

That $r^{\prime}(r,v)$ is determined uniquely in the situation can be seen
explicitly if it is considered that the problem can be restricted to a
$2$-dimensional problem since Euclidean geodesic planes are transferred to
Euclidean geodesic planes in normal coordinates. The geometric PDE (3.1)
restricts to the $2$-dimensional submanifold. Then, by considering the radial mapping

\begin{center}
$f_{r^{\prime}}:%
\begin{pmatrix}
x\\
y
\end{pmatrix}
\rightarrow\frac{r^{\prime}(r,v)}{r}%
\begin{pmatrix}
x\\
y
\end{pmatrix}
=r^{\prime}(r,v)v$,

$r=\sqrt{x^{2}+y^{2}}$, $v=\frac{1}{r}%
\begin{pmatrix}
x\\
y
\end{pmatrix}
$
\end{center}

in the $2$-dimensional situation, the $r^{\prime}(r,v)$ in (3.5) is indeed
uniquely determined by the (infinitesimally volume preserving) global equation

\begin{center}
\hfill$Vol_{\mathbb{R}^{n}}(f_{r^{\prime}}^{-1}(seg_{\bigtriangleup
\varphi,\bigtriangleup r}(B_{r^{\prime}})))=Vol_{g}(seg_{\bigtriangleup
\varphi,\bigtriangleup r}(B_{r^{\prime}}))$\hfill(5.4)
\end{center}

where $B_{r^{\prime}}\subset\mathbb{R}^{n}$ is the Euclidean geodesic ball in
normal coordinates and $seg_{\bigtriangleup\varphi,\bigtriangleup r^{\prime}%
}(B_{r^{\prime}})\subset B_{r^{\prime}}$ denotes a spherical segment with
angle and radius

\begin{center}
$\bigtriangleup\varphi=\varphi_{2}-\varphi_{1}\in\lbrack0,2\pi]$,
$\bigtriangleup r^{\prime}=r_{2}^{\prime}-r_{1}^{\prime}>0$.
\end{center}

(5.4) can be stated simplified with\ polar coordinates

\begin{center}
$f_{polar}^{-1}(r,\varphi)=%
\begin{pmatrix}
r\cos(\varphi)\\
r\sin(\varphi)
\end{pmatrix}
=%
\begin{pmatrix}
x\\
y
\end{pmatrix}
$, $r>0$, $\varphi\in\lbrack0,2\pi)$,
\end{center}

in preimage and image space as follows

\begin{center}
\hfill$\underset{f_{r^{\prime}}^{-1}(seg_{\bigtriangleup\varphi,\bigtriangleup
r}(B_{r^{\prime}}^{g})}{%
{\textstyle\int}
}ldld\phi=\overset{\varphi_{2}}{\underset{\varphi_{1}}{%
{\textstyle\int}
}}\overset{r_{2}^{\prime}}{\underset{r_{1}^{\prime}}{%
{\textstyle\int}
}}l^{\prime}g(l^{\prime},\phi^{\prime})dl^{\prime}d\phi^{\prime}$,\hfill(5.5)
\end{center}

where $g=\sqrt{\det(g_{ij})^{ij}}$ is the $2$-dimensional restriction volume
element in normal coordinates and $\det(Df_{polar}^{-1})=r$. The solution of
(5.5) is uniquely obtained by

\begin{center}
\hfill$\frac{dr^{\prime}}{dr}=\frac{r}{r^{\prime}g}$\hfill(5.6)
\end{center}

because

\begin{center}
$\det(Df_{r^{\prime}}^{-1})=\frac{dr}{dr^{\prime}}\frac{r}{r^{\prime}}%
\overset{\text{(5.6)}}{=}g$,
\end{center}

holds exactly if (5.6) is obeyed, where the first equality can be seen by
elementary analysis for a $2$-dimensional radial mapping with

\begin{center}
$\det(Df_{polar}\circ Df_{r^{\prime}}^{-1}\circ Df_{polar}^{-1})=\det
(D(f_{polar}\circ f_{r^{\prime}}^{-1}\circ f_{polar}^{-1}))=\frac
{dr}{dr^{\prime}}$.
\end{center}

The radial mapping induces the volume element in (5.5).

The integral equation (5.4) determines $r^{\prime}(r,v)$, where the
$2$-dimensional unit vector $\left\Vert v\right\Vert =1$ has polar coordinates
$\varphi$.

Moreover, it follows that in the $2$-dimensional case the image measure of the
Lebesgue measure of synthetical coordinates is the inner volume of the
Riemannian manifold.\hfill\hfill\ \ \ \ \ \ \ \ \ \ $\blacksquare$

Since $\Theta_{g}$ is an isometry it commutes with the musical isomorphisms
$.^{\#}$ respectively $.^{\flat}$ and the application of the diffeomorphism
$D\Theta_{g}$ on a curve equivalence class can be defined as in lemma 2.4, II.

\textbf{Definition 5.3 }(Renormalized curve transportation.)\ Let $\Theta_{g}$
denote the Levy-Civita metrical distortion (without turbulence [1,2,3]), which
is defined for the complement of $cutlocus(p)$ as image. Let

\begin{center}
$c(s)\in$ $\mathbb{R}_{synth}^{n}$
\end{center}

denote a curve with $s\in\lbrack0,s_{\max}]$.

The invertible reparametrization $s\rightarrow t$ such that

\begin{center}
\hfill$\left\langle \frac{d}{dt}\Theta_{g}(c(s(t))),\frac{d}{dt}\Theta
_{g}(c(s(t)))\right\rangle ^{g}\overset{!}{=}1$\hfill(5.7)
\end{center}

is called \emph{renormalized curve transportation}, if $c$ is an Euclidean
unit speed curve in $\mathbb{R}_{synth}^{n}$.\hfill\hfill$\blacksquare$

Because $\Theta_{g}$ is an isometry, the covariant differential has a
contravariant curve transportation version. The Euler-scheme of a curve on the
manifold is "lifted in synthetical coordinates". Every curve segment is
transformed to normal coordinates with origin in the start point of the
segment and time changed to the length parameter of $\mathbb{R}_{synth}^{n}$.

\textbf{Remark 5.4 }(I) (Computation of the radial isometry property for
$\Theta_{g}$.) Let $r$ denote the radial parameter in $\mathbb{R}_{synth}^{n}%
$, $r^{\prime}$ the radial parameter on $T_{p}M$ and $r^{\prime\prime
}=r^{\prime}\circ((\exp_{p}^{-1})^{1},..,(\exp_{p}^{-1})^{n})$ the radial
parameter on the manifold. If $e_{r}$ denotes the radial row $n$-tuple then

\begin{center}
$e_{r}=dr=dr^{\prime}=dr^{\prime\prime}$
\end{center}

on the level of row $n$-tuples according to proposition 2.1.

In order to compute the radial isometry property

\begin{center}
$D\Theta_{g}[r]=[r^{\prime\prime}]$
\end{center}

with the decomposition in lemma 2.4, III, the classical differential of the
radial mapping

\begin{center}
$\exp_{p}^{-1}\circ\Theta_{g}:\mathbb{R}_{synth}^{n}\rightarrow$
$\mathbb{R}^{n}$
\end{center}

is computed by the classical covariant matrix transformation in (2.5) as follows

\begin{center}
$D(\exp_{p}^{-1}\circ\Theta_{g})dr=\frac{dr}{dr^{\prime}}dr^{\prime}$.
\end{center}

With lemma 2.4, I

\begin{center}
$D\exp_{p}[r^{\prime}]=\frac{dt}{ds}\partial r^{\prime}$
\end{center}

and according to theorem 3.1, III

\begin{center}
$\frac{dt}{ds}=(\frac{dr}{dr^{\prime}})^{-1}$.
\end{center}

Such that

\begin{center}
$D\Theta_{g}[r]=D\exp_{p}\circ D(\exp_{p}^{-1}\circ\Theta_{g})=\frac
{dr}{dr^{\prime}}(\frac{dr}{dr^{\prime}})^{-1}\partial r^{\prime}=\partial
r^{\prime}=[r^{\prime\prime}]$.
\end{center}

(II) (Computation of the isometry property of the exterior differential for
radial symmetric manifolds.) Let the manifold $M$ be radial symmetric with
origin $p$, i.e.

\begin{center}
$r^{\prime}(r,v)=r^{\prime}(r)$, $\forall v$.
\end{center}

In order to compute that the exterior differential $D\Theta_{g}$ is an
isometry on an elementary level with definition 2.2, it has to be considered
that $D\Theta_{g}$ is a classical covariant inner transformation inclusive the
reparametrization according to lemma 2.4. If it is possible to show that the
inner Jacobian of $\exp_{p}^{-1}\circ\Theta_{g}$ inclusive the global
reparametrization $s\rightarrow t$ transforms orthonormal frames in
synthetical coordinates to orthonormal frames in Riemannian coordinates, it
follows that $\Theta_{g}$ is an isometry with respect to the exterior differential.

The following is used

\begin{quote}
-The geodesically radial contraction is $\frac{d}{dr}r^{\prime}(r)$.\hfill\ \ \ \ \ \ \ \ \ \ 

-The situation can be simplified to the $2$-dimensional case.

-The metrically complementary contraction is determined by the $\varkappa$ in (3.1).
\end{quote}

The first assertion follows by construction according to I. The second follows
because Euclidean geodesic planes are transferred to geodesic planes in
coordinates with origin $p$. The last assertion follows by the inspection of
special coordinate systems. Geodesic polar coordinates $(r^{\prime},\varphi)$
can be introduced in coordinates and the covectors $(\partial_{r^{\prime}%
},\partial_{\varphi})$ have a covector $2$-tuple representation.
$\partial_{r^{\prime}}$ and $\partial_{\varphi}$ can be computed as covector
$2$-tuples on the $(x,y)$ coordinate system with $f_{polar}:(x,y)\rightarrow
(r^{\prime},\varphi)$ by

\begin{center}
$dr^{\prime}(x,y)=\frac{1}{\left\vert (x,y)\right\vert }(x,y)$
\end{center}

respectively

\begin{center}
$d\varphi(x,y)=(-r^{\prime}\sin(\varphi),r^{\prime}\cos(\varphi))(x,y)$.
\end{center}

In normal coordinates $\exp_{p}^{-1}$ the application of $\varkappa$ to
vectors $v=v_{1}\partial_{r^{\prime}}+v_{2}\partial_{\varphi}$ has the form

\begin{center}
\hfill$\varkappa(v)=\frac{\left\langle \partial_{r^{\prime}},v\right\rangle
_{\mathbb{R}^{2}}}{\left\Vert \partial_{r^{\prime}}\right\Vert _{\mathbb{R}%
^{2}}^{2}}\partial_{r^{\prime}}+\frac{1}{g}\frac{\left\langle \partial
_{\varphi},v\right\rangle _{\mathbb{R}^{2}}}{\left\Vert \partial_{\varphi
}\right\Vert _{\mathbb{R}^{2}}^{2}}\partial_{\varphi}$,\hfill(5.8)
\end{center}

because normal coordinates are radially length preserving, the manifold is
radial symmetric at the outset and $\det(\varkappa)=\frac{1}{g}$, with the
$2$-dimensional restriction volume element $g$ in normal coordinates, such
that the respective eigenvalues can be concluded necessarily by the isometry
property (3.3).

Because in the $2$-dimensional situation in (5.5) the inner differential of
the radial mapping leads necessarily to the same volume form and the parameter
$\varphi$ is not changed

\begin{center}
$f_{polar}\circ\exp_{p}^{-1}\circ\Theta_{g}\circ f_{polar}^{-1}(r,\varphi
)=(r^{\prime},\varphi)$,
\end{center}

it has the Jacobian

\begin{center}
$D(f_{polar}\circ\exp_{p}^{-1}\circ\Theta_{g}\circ f_{polar}^{-1})(r,\varphi)=%
\begin{pmatrix}
\frac{r}{r^{\prime}g} & 0\\
0 & 1
\end{pmatrix}
$,
\end{center}

because $\det(Df_{polar}^{-1})=r$. It follows that the mapping $r^{\prime}(r)$ obeys

\begin{center}
$\frac{dr^{\prime}}{dr}=\frac{r}{r^{\prime}g}$.
\end{center}

$\frac{1}{r}\partial\varphi$ is a unit covector, it is an Euclidean unit
covector in the preimage space of $\Theta_{g}$, that is, synthetical
coordinates. Analogously, $\frac{1}{r^{\prime}}\partial\varphi$ is an
Euclidean unit covector in the normal coordinate $\mathbb{R}^{n}$ and because
of the form of $\varkappa$ in (5.8), the vector $\frac{1}{r^{\prime}g}%
\partial\varphi$ is a unit vector with respect to the Riemannian metric. The
covariant inner Jacobian of the radial mapping transforms the frame
$\partial\varphi(r,\varphi)$ to $\partial\varphi(r^{\prime},\varphi)$ where
the angle $\varphi$ is the same in image and preimage space, i.e., the inner
Jacobian $D(\exp_{p}^{-1}\circ\Theta_{g})$ is $id$ in the $\varphi$-component
for covector tuples. With lemma 2.4, the differential is

\begin{center}
$D\Theta_{g}\frac{1}{r}[\varphi]=\frac{1}{r^{\prime}g}\partial\varphi=\frac
{1}{r^{\prime}g}[\varphi((\exp_{p}^{-1})^{1},..,(\exp_{p}^{-1})^{n})]$,
\end{center}

where the angle parameter is the same, but the covectors refer to different
domains. That is, $D\Theta_{g}$ is an isometry in the $\varphi$-component.
$D\Theta_{g}$ is by construction radially isometric according to I. Since
$\partial_{r^{\prime}}$ and $\partial_{\varphi}$ are metrically orthogonal it
follows that $D\Theta_{g}$ obeys the isometry property.\hfill\hfill
$\blacksquare$

\begin{quote}
\textbf{6 Exterior Geometry of the 2-Sphere}
\end{quote}

A natural example for a metrical distortion as radial symmetric projection is
provided by the embedding of the $2$-dimensional sphere $S^{2}$ in
$\mathbb{R}^{3}$

\begin{center}
\hfill$S^{2}\subset\mathbb{R}^{3}$,\hfill(6.1)

$S^{2}:=\{(x,y,z)^{T}\in\mathbb{R}^{3}|$ $\left\Vert (x,y,z)^{T}\right\Vert
_{\mathbb{R}^{3}}=1\}$. \ \ \ \ 
\end{center}

The double tangential space in the north pole $n:=(0,0,1)^{T}$

\begin{center}
\hfill$\mathbb{R}_{synth}^{n}=T_{0}T_{n}S^{2}$, $T_{n}S^{2}=\{(x,y,z)^{T}%
\in\mathbb{R}^{3}|$ $z=1\}\subset\mathbb{R}^{3}$,\hfill\ \ \ \ \ \ \ \ \ \ 
\end{center}

serves as synthetical coordinate system with $n$ as origin.

A natural inner coordinate system is the orthonormal projection on
$\mathbb{R}^{3}|_{z=0}$.

\textbf{Remark 6.1 }The orthonormal projection

\begin{center}
\hfill$\Pi_{S^{2}}:S^{2}\rightarrow\mathbb{R}^{3}|_{z=0}:=\{(x,y,z)^{T}%
\in\mathbb{R}^{3}|$ $z=0\}$,\hfill(6.2)

$(x,y,z)^{T}\rightarrow(x,y,0)^{T}$ \ \ \ \ \ 
\end{center}

as chart of $S^{2}$ on $(\mathbb{R}^{3}|_{z=0},\left\langle .,.\right\rangle
^{g})$ provides a natural inner coordinate system with image

\begin{center}
$image(\Pi_{S^{2}})=\{(x,y,0)^{T}\in\mathbb{R}^{3}|_{z=0}|$ $\left\Vert
(x,y)^{T}\right\Vert _{\mathbb{R}^{2}}<1\}$
\end{center}

and the open half sphere

\begin{center}
$S^{2,\uparrow}:=\{(x,y,z)^{T}\in\mathbb{R}^{3}|$ $z=\sqrt{1-x^{2}-y^{2}%
},(x,y)^{T}\in image(\Pi_{S^{2}})\}$
\end{center}

as preimage.

The induced metrical tensor of $\left\langle .,.\right\rangle _{\mathbb{R}%
^{3}}$ restricted to the half sphere $S^{2,\uparrow}$ is in orthonormal
projection coordinates with

\begin{center}
$D\Pi_{S^{2}}^{-1}=%
\begin{pmatrix}
1 & 0\\
0 & 1\\
\frac{-x}{\sqrt{1-x^{2}-y^{2}}} & \frac{-y}{\sqrt{1-x^{2}-y^{2}}}%
\end{pmatrix}
$
\end{center}

computed by

\begin{center}
\hfill$(g_{ij})^{ij}(x,y)=D\Pi_{S^{2}}^{-1,T}\circ D\Pi_{S^{2}}^{-1}%
(x,y)$\hfill(6.3)

$=%
\begin{pmatrix}
1+\frac{x^{2}}{1-x^{2}-y^{2}} & \frac{xy}{1-x^{2}-y^{2}}\\
\frac{xy}{1-x^{2}-y^{2}} & 1+\frac{y^{2}}{1-x^{2}-y^{2}}%
\end{pmatrix}
$. \ \ \ \ \ 
\end{center}

The mapping

\begin{center}
$(x,y)^{T}\rightarrow(x,y,\sqrt{1-x^{2}-y^{2}})^{T}$
\end{center}

is classically differentiable on $image(\Pi_{S^{2}})$.\hfill\hfill
\ \ \ \ \ \ \ \ \ \ $\blacksquare$

An elementary computation with the trigonometric functions shows that in polar
coordinates $f_{polar}:(x,y)\rightarrow(\widehat{r},\varphi)$ of the
orthonormal projection coordinate system

\begin{center}
\hfill$(g_{ij})^{ij}(\widehat{r},\varphi)=(Df_{polar}^{-1})^{T}\circ
(g_{ij})^{ij}(f_{polar}^{-1}(\widehat{r},\varphi))\circ Df_{polar}^{-1}$\hfill(6.4)

$=%
\begin{pmatrix}
\frac{1}{1-\widehat{r}^{2}} & 0\\
0 & \widehat{r}^{2}%
\end{pmatrix}
$, $\widehat{r}=\sqrt{x^{2}+y^{2}}$, \ \ \ \ \ 
\end{center}

with

\begin{center}
\hfill$g_{\widehat{r}\widehat{r}}(\widehat{r},\varphi)=\frac{1}{1-\widehat
{r}^{2}}$, $\forall\varphi$,\hfill(6.5)
\end{center}

where Euclidean polar coordinates are defined by

\begin{center}
$f_{polar}^{-1}(\widehat{r},\varphi)=%
\begin{pmatrix}
\widehat{r}\cos(\varphi)\\
\widehat{r}\sin(\varphi)
\end{pmatrix}
=%
\begin{pmatrix}
x\\
y
\end{pmatrix}
$, $\widehat{r}>0$, $\varphi\in\lbrack0,2\pi)$,

$((g_{polar})_{ij})^{ij}=(Df_{polar}^{-1})^{T}\circ Df_{polar}^{-1}=%
\begin{pmatrix}
1 & 0\\
0 & \widehat{r}^{2}%
\end{pmatrix}
$.
\end{center}

\textbf{Example 6.2 }(Riemannian exponential mapping as radially length
preserving radial symmetric projection on the sphere.) The Riemannian
exponential mapping $\exp_{n}$ with image $S^{2}\subset\mathbb{R}^{3}$ is now
determined by the radially length preserving radial symmetric projection

\begin{center}
\hfill$\exp_{n}:B_{r^{\prime}<\pi/2}\subset\mathbb{R}^{2}\rightarrow
S^{2}\overset{\Pi_{S^{2}}}{\rightarrow}image(\Pi_{S^{2}})$\hfill(6.6)

$(r^{\prime},\varphi)\underset{\text{projection}}{\overset{\text{radial}%
}{\rightarrow}}(\widehat{r},\varphi)$, \ \ \ \ \ 
\end{center}

where $r^{\prime}$ denotes the length parameter in normal coordinates and
$B_{r^{\prime}<\pi/2}$ denotes the $2$-ball with radius $\pi/2$.

The explicit function $\widehat{r}(r^{\prime})$ is determined by the endpoint
of any radial unit speed curve $c(t)$ in the coordinate system. The unit speed
parameter $t$ equal to $r^{\prime}$ determines the function. That is

\begin{center}
$\sqrt{\left\langle \frac{d}{dt}c(t),\frac{d}{dt}c(t)\right\rangle ^{g}}%
=\sqrt{g_{\widehat{r}\widehat{r}}(c(t))(\frac{d}{dt}c(t))^{2}}\overset{!}{=}%
1$, $\forall t$,
\end{center}

with

\begin{center}
$\frac{1}{1-\widehat{r}(r^{\prime})^{2}}=g_{\widehat{r}\widehat{r}}%
(\widehat{r}(r^{\prime}))$, $\forall\varphi$.
\end{center}

The latter condition is equivalent to the ODE

\begin{center}
$\frac{d}{dr^{\prime}}\widehat{r}(r^{\prime})=\sqrt{1-\widehat{r}(r^{\prime
})^{2}}$, $\widehat{r}(0)=0$, $\widehat{r}>0$.
\end{center}

The ODE has the classical closed solution

\begin{center}
\hfill$\widehat{r}(r^{\prime})=\sin(r^{\prime})$, $r(\pi/2)=1$,\hfill(6.7)
\end{center}

what is easily seen with $\cos(\alpha)^{2}+\sin(\alpha)^{2}\equiv1$,
$\forall\alpha$.

With the Gau\ss 's\ lemma, theorem 5.1, by the uniqueness of the radial
symmetric length preserving mapping, (6.6) is the exponential mapping with the
north pole $n$ as origin.\hfill\hfill\ \ \ \ \ \ \ \ \ \ $\blacksquare$

\textbf{Remark 6.3 }(I) (Determination of $\varkappa$ for the metrical
distortion of the $2$-sphere.) If $(\widehat{r},\varphi)$ parametrize polar
coordinates of the projection coordinate system, then $\varkappa$ in
projection coordinates is determined by the matrix in the projection
coordinate $\mathbb{R}^{n}$ which is determined by the eigenvectors
$\partial_{\widehat{r}}$ and $\partial_{\varphi}$ in coordinates respectively
their eigenvalues $\sqrt{1-\widehat{r}^{2}}$ and $1$. This is seen by radial
length preservation respectively radial symmetry because the orthonormal
projection of $z$-axis level set circles on the sphere to distance circles in
the projection coordinate system is length preserving along the circles, such
that the eigenvalue $1$ for $\partial_{\varphi}$ can be concluded. Concretely,
as explicit analytic expression, $\varkappa$ is defined by

\begin{center}
\hfill$\varkappa(v)=\sqrt{1-\widehat{r}^{2}}\frac{\left\langle \partial
_{\widehat{r}},v\right\rangle _{\mathbb{R}^{2}}}{\left\Vert \partial
_{\widehat{r}}\right\Vert _{\mathbb{R}^{2}}^{2}}\partial_{\widehat{r}}%
+\frac{\left\langle \partial_{\varphi},v\right\rangle _{\mathbb{R}^{2}}%
}{\left\Vert \partial_{\varphi}\right\Vert _{\mathbb{R}^{2}}^{2}}%
\partial_{\varphi}$,\hfill(6.8)

$v=v_{1}\partial_{\widehat{r}}+v_{2}\partial_{\varphi}$ \ \ \ \ \ 
\end{center}

because the expression is linear and it performs the eigenvalue theory which
is unique.

Alternatively, the determination of $\varkappa$ for the metrical distortion
proceeds according to [1, remark 7] by the numerical solution of the parallel
translation ODE in coordinates along radial geodesics.

An inspection of $\varkappa$ in polar coordinates shows that

\begin{center}
$\varkappa\circ\varkappa^{T}=(g^{ij})_{ij}$
\end{center}

in (6.4) is not the diagonal solution.

(II) (Determination of the contraction factor.) If the metrical tensor is
given in projection coordinates by (6.4) then in normal coordinates

\begin{center}
$g_{ks}=(Df_{r^{\prime}}^{-1})_{k}^{i}(Df_{r^{\prime}}^{-1})_{s}^{j}g_{ij}$
\end{center}

where

\begin{center}
$f_{r^{\prime}}:image(\Pi_{S^{2}})\rightarrow T_{0}image(\Pi_{S^{2}})$

$f_{r^{\prime}}(x)=r^{\prime}(\left\vert x\right\vert )\frac{x}{\left\vert
x\right\vert }$
\end{center}

is the coordinate change from projection coordinates to normal coordinates,
which is a radial mapping. $r^{\prime}\circ\exp_{n}^{-1}$ is the length
parameter of $S^{2}$. An elementary computation shows that

\begin{center}
$g_{normal}=\det(g_{ks})^{ks}=\frac{\sin(r^{\prime})}{r^{\prime}}$.
\end{center}

According to (5.6) the differential slip is

\begin{center}
\hfill$\frac{r}{r^{\prime}g}=\frac{r}{\sin(r^{\prime})}$,\hfill(6.9)
\end{center}

where $r$ denotes the associated length parameter in synthetical coordinates
and $g=g_{normal}$ as global scalar.\hfill\hfill
\ \ \ \ \ \ \ \ \ \ $\blacksquare$

\textbf{Example 6.4 }(Metrical distortion as radially volume preserving radial
symmetric projection on the sphere.) The Levy-Civita metrical distortion in
the situation $S^{2}\subset\mathbb{R}^{3}$ is now determined by the radially
volume preserving radial symmetric projection

\begin{center}
\hfill$\Theta_{S^{2}}:T_{0}T_{n}S^{2}|_{B_{r<r_{\max}}}\rightarrow
S^{2}\overset{\Pi_{S^{2}}}{\rightarrow}image(\Pi_{S^{2}})$\hfill(6.10)

$(r,\varphi)\underset{\text{projection}}{\overset{\text{radial}}{\rightarrow}%
}(\widehat{r},\varphi)$, \ \ \ \ \ \ 
\end{center}

where $r$ denotes the length parameter in $\mathbb{R}_{synth}^{2}=T_{n}S^{2}$
and $.|_{B_{r<r_{\max}}}$ denotes the restriction to the interior of a $2$-ball.

If $\triangle\varphi$ denotes the (normalized) arc length of a radial circle
segment of $S^{1}$ in $\mathbb{R}_{synth}^{2}$, then

\begin{center}
$r\triangle\varphi$
\end{center}

is the Euclidean length of the radial projected circle segment with radius
$r>0$ in $\mathbb{R}_{synth}^{2}$.

Analogously

\begin{center}
$\widehat{r}\triangle\varphi$
\end{center}

is the Riemannian length of the radial symmetric circle segment in the
projective coordinate system, what follows by orthonormal projection of
$S^{2}$. The length of these circle segments in $S^{2}$ and in $image(\Pi
_{S^{2}})$ is preserved by orthonormal projection.

The explicit function $\widehat{r}(r)$ is determined by volume preservation in
polar coordinates of synthetical coordinates and polar coordinates of
projection coordinates\ as in (5.5) with the associated Riemannian volume
element in (6.4)

\begin{center}
$\overset{r+\triangle r}{\underset{r}{%
{\textstyle\int}
}}sds\triangle\varphi\overset{!}{=}\overset{\widehat{r}(r+\triangle
r)}{\underset{\widehat{r}(r)}{%
{\textstyle\int}
}}\frac{s^{\prime}}{\sqrt{1-(s^{\prime})^{2}}}ds^{\prime}\triangle\varphi$.
\end{center}

What leads to

\begin{center}
$\frac{1}{2}\frac{(r+\triangle r)^{2}-r^{2}}{\triangle r}\overset{!}{=}%
\frac{-\sqrt{1-\widehat{r}(r+\triangle r)^{2}}+\sqrt{1-\widehat{r}(r)^{2}}%
}{\triangle r}$.
\end{center}

The application of $\underset{\triangle r\rightarrow0}{\lim}$ is a
differentiation which shows that

\begin{center}
$\frac{d}{dr}(-\sqrt{1-\widehat{r}(r)^{2}})\overset{!}{=}r$.
\end{center}

Integration results in

\begin{center}
$\sqrt{1-\widehat{r}(0)^{2}}-\sqrt{1-\widehat{r}(r)^{2}}\overset{!}{=}\frac
{1}{2}r^{2}$.
\end{center}

because of $\widehat{r}(0)=0$ the solution is

\begin{center}
\hfill$\widehat{r}(r)=r\sqrt{1-\frac{1}{4}r^{2}}$.\hfill(6.11)
\end{center}

Which has the inverse

\begin{center}
\hfill$r(\widehat{r})=\sqrt{2}\sqrt{1-\sqrt{1-r^{2}}}$.\hfill(6.12)
\end{center}

Especially

\begin{center}
$r_{\max}=r(1)=\sqrt{2}$.
\end{center}

In order to check the total volume preservation, the surface volume of the
half sphere $S^{2,\uparrow}$ is

\begin{center}
\hfill$Vol(S^{2,\uparrow})=\frac{1}{2}(4\pi(\widehat{r})^{2})\overset
{\widehat{r}=1}{=}2\pi$.\hfill\ \ \ \ \ \ \ \ \ \ 
\end{center}

The volume of the preimage of $\Theta_{S^{2}}$ is the volume of a $2$-ball
with radius $r_{\max}=\sqrt{2}$

\begin{center}
$Vol(\Theta_{S^{2}}^{-1}(S^{2,\uparrow}))=\pi r^{2}\overset{r=\sqrt{2}}{=}%
2\pi$.
\end{center}

That is, total volume preservation holds true.

As in remark 5.4, $\Theta_{S^{2}}$ is an isometry in the arc angle direction,
for every $r>0$ the Euclidean unit speed circle segment

\begin{center}
$s\rightarrow r%
\begin{pmatrix}
\cos(\varphi(s)+\varphi_{0})\\
\sin(\varphi(s)+\varphi_{0})
\end{pmatrix}
\in\mathbb{R}_{synth}^{2}$,
\end{center}

where

\begin{center}
$\varphi(s)=\frac{s}{r}$, $s\in\mathbb{R}$
\end{center}

is transferred by the renormalized curve transportation in definition 5.3 to
the $g$-unit speed circle segment

\begin{center}
$t\rightarrow\widehat{r}(r)%
\begin{pmatrix}
\cos(\varphi(t)+\varphi_{0})\\
\sin(\varphi(t)+\varphi_{0})
\end{pmatrix}
\in image(\Pi_{S^{2}})$,
\end{center}

where

\begin{center}
$\varphi(t)=\frac{t}{\widehat{r}}$, $t\in\mathbb{R}$.
\end{center}

$\Theta_{S^{2}}$ is an isometry in the orthonormal radial direction by construction.

With the generalized Gau\ss 's\ lemma, theorem 5.2, by the uniqueness of the
radial symmetric volume preserving mapping, $\Theta_{S^{2}}$ is the metrical
distortion with the north pole $n$ as origin.\hfill\hfill
\ \ \ \ \ \ \ \ \ \ $\blacksquare$

\textbf{Remark 6.5 }(I) (Exterior differentiability.) $\Theta_{S^{2}}$ is an
isometry from an Euclidean space to a space with a nonvanishing curvature
tensor, $\Theta_{S^{2}}$ can not be classically differentiable. That is, inner
differentiability cannot cover the problem.

(II) (Almost global existence.) The preimage of the metrical distortion is
restricted by the restricted domain of the orthonormal coordinate projection
in the construction of example 6.4. However, normal coordinates are almost
global definable. If (radial symmetric) normal coordinates would have been
chosen as image domain, the radial symmetric construction of the metrical
distortion could have been continued almost global up to the cut locus, which
is the south pole.

(III) (Lie-group $S^{3}$.) [4, theorem 11] states that the metrical distortion
of the $3$-sphere $S^{3}$, which is a Lie-group, is the matrix exponential
function. This assertion has to consider a differential slip with a revision
of Gau\ss 's lemma as in [4,5]. According to theorem 5.2 and the submanifold
property $S^{2}\subset S^{3}$, the metrical distortion $\Theta_{S^{2}}$ is the
restriction since the isometry property (3.3) necessarily restricts. Since the
renormalized metrical distortion

\begin{center}
$\Theta_{S^{3}}(rv)=\exp_{id}(r^{\prime}(r)v)$,

$\left\Vert v\right\Vert =1$, $v\in T_{id}S^{3}\simeq T_{0}T_{id}S^{3}$
\end{center}

is the exponential function and the integral curve respectively the left
invariance argumentation with the autodiffeomorphic differential of the left multiplication

\begin{center}
$L_{x}:S^{3}\rightarrow S^{3}$, $x\in S^{3}$
\end{center}

refers to the parametrization by $r^{\prime}$ in the manifold, the statement
in the theorem should note the equality of the Riemannian exponential function
and the lie group matrix exponential function. However, for the isometry
property of the mapping $\exp_{id}:T_{id}S^{3}\rightarrow S^{3}$ as in lemma
2.4, I the reparametrization due to the metrical distortion is of vital
interest for the differential topology of Lie groups [4, theorem
15].\hfill\hfill\ \ \ \ \ \ \ \ \ \ $\blacksquare$

\begin{quote}
\textbf{Bibliography}

[1] S. V\"{o}llinger, Metrical Distortions and Geometric Partial Differential
Equations, \emph{JP J.\ Geom. Top. }\textbf{14, }(2013), 61-86.

[2] S. V\"{o}llinger, Turbulence Theory, \emph{JP J.\ Geom. Top. }\textbf{15,
}(2014), 35-64.

[3] S. V\"{o}llinger, Metrical Distortions and Invariant Hydrodynamics,
\emph{Far East J. Appl. Math.} \textbf{94}, (2016), 377-394.

[4] S. V\"{o}llinger, Metrical Distortions and Lie Groups, \emph{JP J.\ Geom.
Top. }\textbf{15, }(2018), 35-64.

[5] M. do Carmo, Riemannian Geometry, Birkh\"{a}user, Basel, (1992).

[6] K. J\"{a}nich, Vektoranalysis, Springer, Berlin (1993), Second\ Edition.

[7] J. Elstrodt, Ma\ss .- und Integrationstheorie, Springer, Berlin (1996).
\end{quote}

\end{document}